\documentclass[12pt,a4paper]{article}

\usepackage[american]{babel}
\usepackage[latin1]{inputenc}
\usepackage{graphics}
\usepackage{epsfig}
\usepackage{amsthm,amsfonts,amsmath}
\usepackage[T1]{fontenc}
\usepackage{ulem}
\usepackage{indentfirst}
\usepackage{color}
\usepackage{makeidx}
\usepackage{amssymb}

\newcommand{\cqd}{\begin{flushright}\hspace{5pt}\rule{5pt}{5pt}\end{flushright}}

\newcommand{\real}{\mathbb{R}}
\newcommand{\h}{\mathbb{H}^2}
\newcommand{\s}{\mathbb{S}^{2}}

\newcommand{\dem}{{\bf Proof.}}

\newtheorem{thm}{Theorem}[section]
\newtheorem{prop}{Proposition}[section]
\newtheorem{cor}{Corollary}[section]
\newtheorem{lem}{Lemma}[section]

\begin{document}

\begin{center} {\LARGE \bf Simons' Type Equation in $\s\times\real$ and $\h\times\real$ and Applications}\end{center}

\begin{center}
\textit{Márcio Batista}\footnote{The author was partially
supported by CNPq, at UFAL, Brazil.}
\end{center}

\begin{abstract}
Equations of Simons type are presented. They are satisfied by a
pair of special operators associated to the immersion
$\Sigma^2\looparrowright M^2(c)\times\real$ with constant mean
curvature. Some immersions are characterized.
\end{abstract}

$\it{Keywords}:$ Constant Mean Curvature; Simons' Type Equation.

$2000 \ \  Mathematics \ \ Subject \ \ Classification$: 53E40.

\section{Introduction}
In 1994, using the traceless second fundamental form $\phi = A -
HI$ associated to an immersed hypersurface
$M^n\looparrowright\mathbb{S}^{n+1}$, H. Alencar and M. do Carmo,
see $[AdC]$, proved that \vspace{0.5cm}

$\bf{Theorem.}$ Let $M^n\looparrowright\mathbb{S}^{n+1}$ be an
immersed hypersurface. If $M^n$ is compact and orientable with
constant mean curvature $H$ and
\begin{equation*}
|\phi|^2 \leq B_H,
\end{equation*}
where $B_H$ is the square of the positive root of
\begin{equation*}
P_H(x)=x^2+\dfrac{n(n-2)}{\sqrt{n(n-1)}}Hx-n(H^2+1).
\end{equation*}
Then:
\begin{enumerate}
\item[(a)] Either $|\phi|^2 = 0$ (and $\Sigma^2$ is totally umbilic) or $|\phi|^2 =
B_H$.
\item[(b)] The $H(r)$-tori $\mathbb{S}^{n-1}(r) \times \mathbb{S}^1 ( 1
- r^2)$ with $r^2 \leq \dfrac{n-1}{n}$ are the only hypersurfaces
with constant mean curvature H and $|\phi|^2 = B_H$.
\end{enumerate}
\vspace{0.5cm}

Motivated by this result we decided to study this problem for
surfaces in $M^2(c)\times\real$ with $c=\pm1$, where $M^2(-1)=\h$
and $M^2(1)=\s$.

We begin by using the traceless second fundamental form $\phi$
associate to an immersed surface $\Sigma^2\looparrowright
M^2(c)\times\real$.

Then we study, in section 3, a special tensor $S$ defined by
\begin{equation}\label{s}
SX= 2HAX-c<X,T>T+\dfrac{c}{2}(1-\nu^2)X-2H^2X,
\end{equation}
where $X \in T_p\Sigma$, $A$ is the Weingarten operator associated
to the second fundamental form, $H$ is the mean curvature, $T$ is
the tangential component of the parallel field $\partial_t$, tangent
to $\real$ in $M^2(c)\times\real$, and $\nu = <N,\partial_t>$.

This operator satisfies Codazzi's equation, provided $H$ is constant,
and its trace vanishes. We remark that any surface with $|S|=0$ and
constant mean curvature is very interesting, because the $(2,0)$-part
of the quadratic differential Q,
\begin{equation*}
Q(X,Y)=2H<AX,Y>-c<X,\partial_t><Y,\partial_t>,
\end{equation*}
of these surfaces is null. In $[AR]$ section 2, Abresch and
Rosenberg described four distinct classes of complete, possibly
immersed, constant mean curvature surfaces
$\Sigma^2\looparrowright M^2(c)\times\real$ with vanishing
(2,0)-part of quadratic differential Q.

Three of these classes are comprised by the embedded rotationally
invariant surfaces; $(i)$ they are the constant mean curvature
spheres $S^2_H \looparrowright M^2(c)\times\real$ of Hsiang and
Pedrosa, see $[Hs]$ section 3 to the case $\h\times\real$ and $[PR]$
to the case $\s\times \real$, $(ii)$ their non-compact cousins
$D^2_H$, and $(iii)$ the surfaces $C^2_H$ of catenoidal type. The
fourth class is comprised of certain orbits $P^2_H$ of the
2-dimensional solvable groups of isometries of $M^2(c)\times \real$.
In this paper we will call these surface described in $[AR]$
Abresch-Rosenberg surfaces.
\\
\\
Then we obtain an equation of Simons type to $S$ and apply it in some
particular cases:

\begin{thm}
Let $\Sigma^2 \looparrowright M^2(c)\times\real$ be an immersed
surface with non zero constant mean curvature $H$ and $S$ as
defined in $(\ref{s})$. Then,
\begin{equation*}
\dfrac{1}{2}\Delta|S|^2=|\nabla
S|^2-|S|^4+|S|^2\left(\dfrac{5c\nu^2}{2}-\dfrac{c}{2}+2H^2-\dfrac{c}{H}<ST,T>\right)
\end{equation*}
\begin{equation*}
+c|ST|^2-\dfrac{1}{4H^2}<ST,T>^2.
\end{equation*}
\end{thm}

\vspace{0.5cm}

 Let us consider the polynomial
$p_H(t)=-t^2-\dfrac{1}{H}t+\left(\dfrac{4H^2-1}{2}\right)$.  When
$H$ is greater than one half then there is a positive root for
$p_H$. Let $L_H$ be such the positive root. One has:
\begin{thm}
Let $\Sigma^2 \looparrowright \s\times\real$ be a complete
immersed surface with constant mean curvature $H$ greater than one
half. Assume that
\begin{equation*}
\sup_{\Sigma}|S| < L_H.
\end{equation*}
Then $\Sigma^2 = S_H^2$, i.e, $\Sigma^2$ is a Hsiang-Pedrosa
sphere.
\end{thm}
Remark.1: The number $L_H$ is
$\dfrac{4H^2-1}{\sqrt{8H^4-2H^2+1}-1}$.
\vspace{0.25cm}

Let us consider the polynomial
$q_H(t)=-t^2-\dfrac{1}{H}t+\left(\dfrac{8H^4-12H^2-1}{4H^2}\right)$.
When $H$ greater than $\sqrt{\dfrac{12+\sqrt{176}}{16}}$, then
there is a positive root for $q_H$. Let $M_H$ be such the positive root.

\begin{thm}
Let $\Sigma^2 \looparrowright \h\times\real$ be a complete
immersed surface with constant mean curvature $H$ greater than
$\sqrt{\dfrac{12+\sqrt{176}}{16}}\approx 1.25664$. Assume that
\begin{equation*}
\sup_{\Sigma}|S| < M_H.
\end{equation*}
Then $\Sigma^2$ is an Abresch-Rosenberg surface.
\end{thm}

Remark.2: The number $M_H$ is
$\dfrac{8H^4-12H^2-1}{2H(\sqrt{8H^4-12H^2-1}+1)}$.

\vspace{0.25cm}

Remark.3: We have obtained some results in the compact case in section 4.

\vspace{0.25cm}


$\bf{Acknowledgements}$.  This is part of my Doctoral Dissertation
at IMPA. I want to thank M. do Carmo for his orientation and H.
Alencar for helpful suggestions and his co-orientation.


\section{Preliminaries}

Let $\Sigma^2 \looparrowright M^3$ be an immersed surface. Let
$\overline{\nabla}$ denote the Levi-Civitá connection on $M^3$ and
let $\nabla$ denote the Levi-Civitá connection on $\Sigma$ for the
induced metric.

Generally speaking, objects defined on $M^3$ will be denoted by
the same symbols as the corresponding objects defined on $\Sigma$
plus a bar over the symbol.

The Riemannian metric extends to natural inner products on space
of tensors and the above connections induce natural covariant
derivatives on tensor fields. For example, for $\{e_1,e_2\}$ a
geodesic frame in $p \in \Sigma^2$ and a tensor $\psi$ on
$\Sigma^2$, we have
\begin{equation*}
\nabla^2\psi(p)=\sum_{i=1}^2(\nabla_{e_i}\nabla_{e_i}\psi)(p).
\end{equation*}
For more details about covariant derivatives on tensor fields see
[S], sections 1 and 2.

 We adopt the following convention for curvature tensor: if $x,y,z \in T_p\Sigma$, we define $R_{x,y}z$ by
\begin{equation*}
R_{x,y}z = R(X,Y)Z(p)=(\nabla_X\nabla_YZ-\nabla_Y\nabla_XZ -
\nabla_{[X,Y]}Z)(p),
\end{equation*}
for any local vector fields which extend the given vectors
$x,y,z$.

The second fundamental form is defined by $\alpha(X,Y) =
(\overline{\nabla}_XY)^{\perp}$ and the Weingarten operator
associated is given by $Av =-(\overline{\nabla}_vN)^{T}$, where $N$
is an unit normal field on $\Sigma^2$. We use the Weingarten operator to
define the operators
\begin{equation*}
\hspace{-1cm}<\overline{R}(A)x,y> :=
\sum_{i=1}^2(-<Ax,\overline{R}_{e_i,y}e_i>-<Ay,\overline{R}_{e_i,x}e_i>
\end{equation*}
\begin{equation}\label{1}
\hspace{2cm}
+<Ay,x><N,\overline{R}_{e_i,N}e_i>-2<Ae_i,\overline{R}_{e_i,x}y> )
\end{equation}
and
\begin{equation*}
<\overline{R}'x,y>:=
\sum_{i=1}^2\{<(\overline{\nabla}_x\overline{R})_{e_i,y}e_i,N>+<(\overline{\nabla}_{e_i}\overline{R})_{e_i,x}y,N>\},
\end{equation*}
where $\{e_1,e_2\}$ is a orthonormal basis of $T_p\Sigma$.
\\

With this notation we have the following results:
\begin{thm}\label{t1}
Let $\Sigma^2 \looparrowright M^3$ be an immersed surface with
constant mean curvature $H$. For any $x,y \in T_p\Sigma$ we have
\begin{equation*}
\hspace{-1cm}<(\nabla^2A)x,y>= -|A|^2<Ax,y> + <\overline{R}(A)x,y>
\end{equation*}
\begin{equation}\label{2}
\hspace{2cm}+<\overline{R}'x,y> +
2H<\overline{R}_{N,x}y,N>+2H<Ax,Ay>.
\end{equation}
\end{thm}
\dem \ \ See Theorem 2 in [B] assumed that co-dimensional is 1. \cqd

We will use the Omori-Yau Maximum Principle. The proof this result
can be found in $[Y]$, Theorem 1.

\begin{thm}[Omori-Yau Maximum Principle]\label{OY}
Let $M$  be a complete Riemannian manifold with Ricci curvature
bounded from below. If $u \in C^{\infty}(M)$ is bounded from
above, then there exist a sequence of points $\{p_j\} \in M$ such
that
\begin{equation*}
\lim_{j\rightarrow\infty}u(p_j)=\sup_M u, \ \ ||\nabla u||(p_j)<
\dfrac{1}{j}, \ and \ \Delta u(p_j) < \dfrac{1}{j}.
\end{equation*}
\end{thm}

\vspace{1cm}

 Let us recall the Gauss' equation to $\Sigma^2$ in
$M^2(c)\times\real$:
\begin{equation*}
R(Y,X)Z=<AX,Z>AY-<AY,Z>AX+c(<X,Z>Y-<Y,Z>X
\end{equation*}
\begin{equation*}
-<Y,T><X,Z>T-<X,T><Z,T>Y
\end{equation*}
\begin{equation}\label{gs}
+<X,T><Y,Z>T+<Y,T><Z,T>X),
\end{equation}
where $X,Y,Z$ in $T_p\Sigma$, $N$ is a unitary normal field on
$\Sigma^2$ and $T$ is the tangential component of the parallel
field $\partial_t$. For more details see $[D]$.

\section{Simons' equation in $M^2(c)\times\real$}

In this section we will obtain an equation of Simons type for the
traceless second fundamental form $\phi$ and for $S$ defined in
$(\ref{s})$.

Let $M^2(c)\times\real$, where $M^2(-1) = \h$ and $M^2(1)= \s$. In
this case we have that $\overline{R}'$=0, because
$M^2(c)\times\real$ is locally symmetric.

In Lemmas $\ref{l1}$ and $\ref{l2}$ we will consider an immersed
surface $\Sigma^2\looparrowright M^2(c)\times\real$ with constant
mean curvature $H$ where $A$ is the Weingarten operator associated
to the second fundamental form on $\Sigma^2$.

\begin{lem}\label{l1}
Denoting the identity by I, we have that
\begin{equation*}\overline{R}(A)=c(5\nu^2-1)A-4cH\nu^2I.\end{equation*}
\end{lem}
\dem \ \ Consider a geodesic frame $\{e_1,e_2\}$ in $p \in
\Sigma^2$ and $T =
\partial_t - \nu N$, where $\nu=<N,\partial_t>$. For $x,y \in T_p\Sigma^2$, we compute
\begin{equation*}
\hspace{-0.5cm}<Ax,\overline{R}_{e_i,y}e_i> =-
c(<y,Ax>-<y,e_i><Ax,e_i>
\end{equation*}
\begin{equation*}
- <y,T><Ax,T>-<e_i,T>^2<y,Ax>
\end{equation*}
\begin{equation*}
+ <e_i,T><Ax,T><e_i,y>+ <y,T><e_i,T><e_i,Ax>).
\end{equation*}
Therefore
\begin{equation*}
\sum_{i=1}^2<Ax,\overline{R}_{e_i,y}e_i> =-
c(2<y,Ax>-\sum_{i=1}^2<y,e_i><Ax,e_i>
\end{equation*}
\begin{equation*}
- 2<y,T><Ax,T>-\sum_{i=1}^2<e_i,T>^2<y,Ax>
\end{equation*}
\begin{equation*}
+ <Ax,T>\sum_{i=1}^2<e_i,T><e_i,y>+
<y,T>\sum_{i=1}^2<e_i,T><e_i,Ax>),
\end{equation*}
i.e,
\begin{equation*}
\sum_{i=1}^2<Ax,\overline{R}_{e_i,y}e_i> = -c(2<y,Ax>-<Ax,y>-
2<y,T><Ax,T>
\end{equation*}
\begin{equation*}
-|T|^2<y,Ax>+ <Ax,T><T,y>+ <y,T><Ax,T>)
\end{equation*}
\begin{equation*}
=-c(<Ax,y>-|T|^2<Ax,y>)=-c(1-|T|^2)<Ax,y>.
\end{equation*}

Thus,
\begin{equation*}\sum_{i=1}^2<Ax,\overline{R}_{e_i,y}e_i>=-c\nu^2<Ax,y>.\end{equation*}
Now
\begin{equation}\label{esp}
<\overline{R}_{e_i,N}e_i,N> =-
c\{<e_i^*,e_i^*><N^*,N^*>-<e_i^*,N^*>^2\},
\end{equation}
where $v^*=v - <v,\partial_t>\partial_t$ for any $v \in
T_p(M^2(c)\times\real)$.
\\
We compute
\\
\begin{equation*}
<e_i^*,e_i^*>=<e_i - <e_i,\partial_t>\partial_t,e_i -
<e_i,\partial_t>\partial_t>
\end{equation*}
\begin{equation*}
=<e_i,e_i>-2<e_i,\partial_t><e_i,\partial_t>+<e_i,\partial_t>^2
\end{equation*}
\begin{equation}\label{v1}
=1-<e_i,\partial_t>^2,
\end{equation}
\begin{equation*}
<N^*,N^*>=<N-\nu\partial_t,N -\nu\partial_t>
\end{equation*}
\begin{equation*}
=<N,N>-2\nu<N,\partial_t>+\nu^2<\partial_t,\partial_t>
\end{equation*}
\begin{equation}\label{v2}
=1-\nu^2,
\end{equation}
\begin{equation*}
<N^*,e_i^*>=<N-\nu\partial_t,e_i - <e_i,\partial_t>\partial_t>
\end{equation*}
\begin{equation*}
<N,e_i>-<e_i,\partial_t><N,\partial_t>-\nu<\partial_t,e_i>+\nu<e_i,\partial_t><\partial_t,\partial_t>
\end{equation*}
\begin{equation} \label{v3}
=\nu<e_i,\partial_t>.
\end{equation}

 Taking $(\ref{v1})$, $(\ref{v2})$ and $(\ref{v3})$ into $(\ref{esp})$ obtain
\begin{equation*}
<\overline{R}_{e_i,N}e_i,N>
=-c\{(1-<e_i,\partial_t>^2)(1-\nu^2)-\nu^2<e_i,\partial_t>^2\}
\end{equation*}

\begin{equation*}
=-c\{1-<e_i,\partial_t>^2-\nu^2+<e_i,\partial_t>^2\nu^2-<e_i,\partial_t>^2\nu^2\}
\end{equation*}
\begin{equation*}
=-c\{1-\nu^2-<e_i,\partial_t>^2\}.
\end{equation*}

 Therefore,
\begin{equation*}
\sum_{i_1}^2<\overline{R}_{e_i,N}e_i,N>
=-c\{2-2\nu^2-\sum_{i_1}^2<e_i,\partial_t>^2\}
\end{equation*}
\begin{equation*}
=-c(2-2\nu^2-|T|^2)=-c(2-2\nu^2-1+\nu^2)=-c(1-\nu^2).
\end{equation*}
Next, we observe that
\begin{equation*}
\hspace{-0.5cm}<Ae_i,\overline{R}_{e_i,x}y> =
-c\{<e_i,y><Ax,e_i>-<y,x><Ae_i,e_i>
\end{equation*}
\begin{equation*}
-<x,T><Ae_i,T><e_i,y>-<e_i,T><y,T><e_i,Ax>
\end{equation*}
\begin{equation*}
+ <e_i,T><Ae_i,T><x,y>+ <y,T><x,T><e_i,Ae_i>\},
\end{equation*}
i.e,
\begin{equation*}
\sum_{i=1}^2<Ae_i,\overline{R}_{e_i,x}y> =
-c\{\sum_{i=1}^2<e_i,y><Ax,e_i>-<y,x>\sum_{i=1}^2<Ae_i,e_i>
\end{equation*}
\begin{equation*}
-<x,T>\sum_{i=1}^2<Ae_i,T><e_i,y>-<y,T>\sum_{i=1}^2<e_i,T><e_i,Ax>
\end{equation*}
\begin{equation*}
+ <x,y>\sum_{i=1}^2<e_i,T><Ae_i,T>+
<y,T><x,T>\sum_{i=1}^2<e_i,Ae_i>\}.
\end{equation*}
Therefore,
\begin{equation*}
\sum_{i=1}^2<Ae_i,\overline{R}_{e_i,x}y>=-c\{<Ax,y>-2H<x,y>-<Ay,T><x,T>
\end{equation*}
\begin{equation*}
-<Ax,T><y,T>+<AT,T><x,y>+2H<x,T><y,T>\}.
\end{equation*}
Using expression $(\ref{1})$ we have that
\begin{equation*}
<\overline{R}(A)x,y>=c\nu^2<Ax,y>+c\nu^2<x,Ay>-<Ax,y>c(1-\nu^2)+2c\{<Ax,y>
\end{equation*}
\begin{equation*}
-2H<x,y>-<Ay,T><x,T>-<Ax,T><y,T>+<AT,T><x,y>
\end{equation*}
\begin{equation*}
+2H<x,T><y,T>\},
\end{equation*}
i.e,
\begin{equation*}
<\overline{R}(A)x,y>=3c\nu^2<Ax,y>+c<Ax,y>-2c<x,T><Ay,T>
\end{equation*}
\begin{equation*}
-2c<Ax,T><y,T>+2c<AT,T><x,y>-4cH<x,y>
\end{equation*}
\begin{equation}\label{esp2}
+4cH<x,T><y,T>.
\end{equation}
\\
If $T=0$ at all points, then we have a slice
$M^2(c)\times\{t_0\}$. Assume that $T\neq 0$ and choose an oriented positive
orthogonal frame $\{T, E \}$, where $T$ is non zero. Then evaluating
$(\ref{esp2})$ in the frame $\{T,E\}$ we obtain

\begin{equation*}<\overline{R}(A)T,T>=3c\nu^2<AT,T>+c<AT,T>-2c<AT,T>|T|^2-4cH|T|^2(1-|T|^2)\end{equation*}
\begin{equation*}=3c\nu^2<AT,T>+c<AT,T>-2c<AT,T>(1-\nu^2)-4cH|T|^2\nu^2\end{equation*}
\begin{equation*}=5c\nu^2<AT,T>-c<AT,T>-4cH\nu^2<T,T>,\end{equation*}

\begin{equation*}<\overline{R}(A)T,E>=3c\nu^2<AT,E>+c<AT,E>-2c<AT,E>\end{equation*}
\begin{equation*}=5c\nu^2<AT,T>-c<AT,T>=<\overline{R}(A)E,T>,\end{equation*}
and
\begin{equation*}<\overline{R}(A)E,E>=3c\nu^2<AE,E>+c<AE,E>+2c<AT,T>|E|^2-4cH|E|^2\end{equation*}
\begin{equation*}=3c\nu^2<AE,E>+c<AE,E>+2c(2H|E|^2-<AE,E>)|E|^2-4cH|E|^2\end{equation*}
\begin{equation*}=3c\nu^2<AE,E>+c<AE,E>-4cH|E|^2(1-|E|^2)-2c<AE,E>|E|^2\end{equation*}
\begin{equation*}=3c\nu^2<AE,E>+c<AE,E>-2c<AE,E>(1-\nu^2)-4cH\nu^2<E,E>\end{equation*}
\begin{equation*}=5c\nu^2<AE,E>-c<AE,E>-4cH\nu^2<E,E>.\end{equation*}
\newpage
It follows that
\begin{equation*}
\overline{R}(A)=c(5\nu^2-1)A-4cH\nu^2I,
\end{equation*}
where $T$ is non zero. Using continuity it can be extended to points
where $T=0$. \cqd
\begin{lem}\label{l2}
$<\overline{R}_{N,x}y,N>=-c\{<x,T><y,T>-<x,y><T,T>\}.$
\end{lem}
\dem \ \ We observe that \\
\begin{equation*}<x^*,y^*>=<x,y>-<x,T><y,T>,\end{equation*}

\begin{equation*}<x^*,N^*>=\nu<x,T>\end{equation*}
and
\begin{equation*}<N^*,N^*>=1-\nu^2,\end{equation*}
where we have used $v^*=v-<v,\partial_t>\partial_t$ for any $v\in
T_p(M^2(c)\times\real)$.
\\
It follows that
\begin{equation*}
<\overline{R}_{N,x}y,N>=-c\{<N^*,x^*><N^*,y^*>-<N^*,N^*><x^*,y^*>\}
\end{equation*}
\begin{equation*}
=-c\{(\nu<x,T>)(\nu<y,T>)-(<x,y>-<x,T><y,T>)<T,T>\}
\end{equation*}
\begin{equation*}
=-c\{\nu^2<x,T><y,T>- <x,y><T,T>+<x,T><y,T>
\end{equation*}
\begin{equation*}
-\nu^2<x,T><y,T>\}=-c\{<x,T><y,T>- <x,y><T,T>\}.
\end{equation*}
This concludes the proof. \cqd

\begin{cor}\label{c1}
Let $\Sigma^2 \looparrowright M^2(c)\times\real$ be an immersed
surface with constant mean curvature $H$ where $A$ is the Weingarten
operator associated to the second fundamental form on $\Sigma^2$.
Then,
\begin{equation*}
<(\nabla^2A)x,y>=-|A|^2<Ax,y>+c(5\nu^2-1)<Ax,y>-4cH\nu^2<x,y>
\end{equation*}
\begin{equation*}
-2cH\{<x,T><y,T>-<x,y><T,T>\}+2H<Ax,Ay>.
\end{equation*}
\end{cor}
\dem \ \ Consider equation $(\ref{2})$
\begin{equation*}
\hspace{-1cm}<(\nabla^2A)x,y>= -|A|^2<Ax,y> + <\overline{R}(A)x,y>
\end{equation*}
\begin{equation*}
\hspace{2cm}+<\overline{R}'x,y> +
2H<\overline{R}_{N,x}y,N>+2H<Ax,Ay>.
\end{equation*}

Now, we use  Lemmas $\ref{l1}$ and $\ref{l2}$ to obtain
\begin{equation*}
\hspace{-1cm}<(\nabla^2A)x,y>= -|A|^2<Ax,y> +
c(5\nu^2-1)<Ax,y>-4cH\nu^2<x,y>
\end{equation*}
\begin{equation*}
-2Hc\{<x,T><y,T>- <x,y><T,T>\}+2H<Ax,Ay>.
\end{equation*}
This concludes the proof. \cqd

Consider two tensors $S,W$ on $\Sigma^2$. We define the inner
product $<S,W>$ in $p \in \Sigma^2$ as
\begin{equation*}
<S,W>=\sum_{i=1}^2<Se_i,We_i>,
\end{equation*}
where $\{e_1,e_2\}$ is an orthonormal basis for $T_p\Sigma$.

\begin{lem}\label{l3}
Let $\Sigma^2\looparrowright M^2(c)\times\real$ be an immersed
surface with constant mean curvature where $A$ is the Weingarten
operator associated to the second fundamental form on $\Sigma^2$.
Then,
\begin{enumerate}
    \item[$(a)$] $<\nabla^2A,I>=0.$
    \item[$(b)$] $<\nabla^2A,A>=-|A|^4+c(5\nu^2-1)|A|^2-8cH^2\nu^2-2cH<AT,T>+4cH^2|T|^2+2Htr(A^3).$
\end{enumerate}
\end{lem}
\dem \ \ Consider $\{e_1,e_2\}$ an orthonormal basis of
$T_p\Sigma$. We use the definition of the inner product between
tensors and the expression in Corollary $\ref{c1}$ to obtain
\begin{equation*}
<\nabla^2A,A>=\sum_{i=1}^2<(\nabla^2A)e_i,Ae_i>=-|A|^2\sum_{i=1}^2<Ae_i,Ae_i>+
\end{equation*}
\begin{equation*}
c(5\nu^2-1)\sum_{i=1}^2<Ae_i,Ae_i>-4cH\nu^2\sum_{i=1}^2<Ae_i,e_i>-2cH\{\sum_{i=1}^2<AT,e_i><e_i,T>
\end{equation*}
\begin{equation*}
-<T,T>\sum_{i=1}^2<Ae_i,e_i>\}+2H\sum_{i=1}^2<A^2e_i,Ae_i>.
\end{equation*}

Therefore,
\begin{equation*}<\nabla^2A,A>=-|A|^4+c(5\nu^2-1)|A|^2-8cH^2\nu^2-2cH<AT,T>+4cH^2|T|^2+2Htr(A^3).\end{equation*}

Using the definition of the inner product and Corollary $\ref{c1}$
we obtain
\begin{equation*}
<\nabla^2A,I>=\sum_{i=1}^2<(\nabla^2A)e_i,e_i>=-|A|^2\sum_{i=1}^2<Ae_i,e_i>+
\end{equation*}
\begin{equation*}
c(5\nu^2-1)\sum_{i=1}^2<Ae_i,e_i>-8cH\nu^2-2cH\{\sum_{i=1}^2<T,e_i><e_i,T>
\end{equation*}
\begin{equation*}
-2<T,T>\}+2H\sum_{i=1}^2<A^2e_i,e_i>.
\end{equation*}
Therefore,
\begin{equation*}<\nabla^2A,I>= -2H|A|^2+c(5\nu^2-1)2H-8cH\nu^2+2cH<T,T>+2H|A|^2\end{equation*}
\begin{equation*}= -2H|A|^2+c(5\nu^2-1)2H-8cH\nu^2+2cH(1-\nu^2)+2H|A|^2\end{equation*}
\begin{equation*}= 10cH\nu^2-2cH-8cH\nu^2+2cH-2cH\nu^2=0,\end{equation*}
where we have used that $\nu^2 + |T|^2=1$. \cqd

\begin{prop}\label{p1}
Let $\Sigma\looparrowright M^2(c)\times\real$ be an immersed
surface with constant mean curvature $H$ and $\phi$ the traceless
second fundamental form, then
\begin{itemize}
    \item[$(a)$] $|\phi|^2 = |A|^2-2H^2.$
    \item[$(b)$] $\nabla\phi = \nabla A.$
    \item[$(c)$] $trA^3 = 3H|\phi|^2+2H^3$.
\end{itemize}
\end{prop}
\dem \ The proof of the item (a) is:
\begin{equation*}
|\phi|^2=<\phi,\phi>=<A-HI,A-HI>=<A,A>-2H<A,I>+H^2<I,I>
\end{equation*}
\begin{equation*}
=|A|^2-4H^2+2H^2=|A|^2-2H^2,
\end{equation*}
where $<A,I>=2H$ and $<I,I>=2$.
\\
To prove the item (b) we consider tangent fields $X,Y$. Then,
\begin{equation*}
(\nabla_X\phi)Y=(\nabla_X A)Y-(\nabla_X(HI))Y=(\nabla_X
A)Y-\nabla_XHI(Y)+H\nabla_XY
\end{equation*}
\begin{equation*}
=(\nabla_X A)Y-H\nabla_XY-X(H)Y+H\nabla_XY=(\nabla_X A)Y,
\end{equation*}
because $H$ is constant.
 \\
Next, the proof to the item(c) is:
\begin{equation*}
tr(A^3)=\sum_{i=1}^2<A^3e_i,e_i>=\sum_{i=1}^2<(\phi+HI)^3e_i,e_i>
\end{equation*}
\begin{equation*}
=\sum_{i=1}^2<(\phi^3+3H\phi^2+3H^2\phi+H^3I)e_i,e_i>=3H|\phi|^2+2H^3,
\end{equation*}
because $tr\phi=tr\phi^3=0$.

 \cqd

\begin{cor}\label{c2}
Let $\Sigma\looparrowright M^2(c)\times\real$ be an immersed
surface with constant mean curvature $H$ and $\phi$ the traceless
second fundamental form. Then
\begin{equation*}
 <\nabla^2\phi,\phi>=
-|\phi|^4+(2H^2+5c\nu^2-c)|\phi|^2-2cH<\phi T,T>.
\end{equation*}
\end{cor}
\dem \ \ We use the  Proposition $\ref{p1}$ and obtain
\begin{equation*}
<\nabla^2\phi,\phi>=<\nabla^2 A,A-HI>=<\nabla^2A,A>-H<\nabla^2
A,I>.
\end{equation*}
Now, we use the Lemma $\ref{l3}$ and obtain
\begin{equation*}
<\nabla^2\phi,\phi>=-|A|^4+c(5\nu^2-1)|A|^2-8cH^2\nu^2+2cH<AT,T>
\end{equation*}
\begin{equation*}
+4cH^2|T|^2+2Htr(A^3).
\end{equation*}

Therefore
\begin{equation*}
<\nabla^2\phi,\phi>=-(|\phi|^2+2H^2)^2+c(5\nu^2-1)(|\phi|^2+2H^2)-8cH^2\nu^2
\end{equation*}
\begin{equation*}
-2cH<(\phi+HI)T,T>+4cH^2|T|^2+2H(3H|\phi|^2+2H^3),
\end{equation*}
which brings us to
\begin{equation*}
<\nabla^2\phi,\phi>=-|\phi|^4-4H^2|\phi|^2-4H^4+c(5\nu^2-1)|\phi|^2+2c(5\nu^2-1)H^2-8cH^2\nu^2
\end{equation*}
\begin{equation*}
-2cH<\phi T,T>-2cH^2|T|^2+4cH^2|T|^2+6H^2|\phi|^2+4H^4.
\end{equation*}
Hence,
\begin{equation*}
<\nabla^2\phi,\phi>=-|\phi|^4+2H^2|\phi|^2+c(5\nu^2-1)|\phi|^2-2cH<\phi
T,T>.
\end{equation*}

\cqd

Now we shall obtain an equation of Simons' type for the traceless
second fundamental form $\phi$:

\begin{thm}\label{t31}
Let $\Sigma^2 \looparrowright M^2(c)\times\real$ be an immersed
surface with constant mean curvature H and $\phi$ is traceless
second fundamental form. Then
\begin{equation*}
\dfrac{1}{2}\Delta|\phi|^2= |\nabla\phi|^2
-|\phi|^4+(2H^2+5c\nu^2-c)|\phi|^2-2cH<\phi T,T>.
\end{equation*}
\end{thm}
\dem \ \ We use that $\dfrac{1}{2}\Delta|\phi|^2=
|\nabla\phi|^2+<\nabla^2\phi,\phi>$ and Corollary $\ref{c2}$. \cqd

Now we evaluate the Laplacian of $|S|^2$ where
\begin{equation}\label{e0}
S = 2HA-c<T,\cdot> T+\dfrac{c}{2}(1-\nu^2)I-2H^2I.
\end{equation}
\begin{thm}
Let $\Sigma^2\looparrowright M^2(c)\times\real$ be an immersed
surface with constant mean curvature and $S$ defined in
$(\ref{e0})$.Then
\begin{equation*}
(\nabla_XS)Y=(\nabla_YS)X,
\end{equation*}
for all tangent fields $X,Y$ on $\Sigma^2$.
\end{thm}
\dem \ \ We consider $(u,v)$ an isothermal parameter to a surface
$\Sigma^2$. Now, we consider the complex parametric to
riemannian metric, $z=u+iv$. Let us set
\begin{equation*}
T_S(X,Y):= (\nabla_XS)Y-(\nabla_YS)X= \nabla_X(SY)-\nabla_Y(SX)-S[X,Y].
\end{equation*}
We will prove that $T_S$ is null. For this, consider

\begin{equation*}
\partial_z=\dfrac{1}{2}(\partial_u-i\partial_v) \ and \
\partial_{\overline{z}}=\dfrac{1}{2}(\partial_u+i\partial_v).
\end{equation*}
Compute,
\begin{equation*}
<T_S(\partial_z,\partial_{\overline{z}}),\partial_z>=
\partial_z<S\partial_{\overline{z}},\partial_z>-<S\partial_{\overline{z}},\nabla_{\partial_z}\partial_z>
\end{equation*}
\begin{equation*}
-\partial_{\overline{z}}<S\partial_z,\partial_z>+ <S\partial_z,\nabla_{\partial_{\overline{z}}}\partial_z>
\end{equation*}
\begin{equation*}
=-Q^{(2,0)}_{\overline{z}}=0,
\end{equation*}
because $Q^{(2,0)}$ is holomorphic and we used that $\nabla_{\partial_z}\partial_{\overline{z}}=0$, $\nabla_{\partial_z}\partial_z=\dfrac{\lambda_{\overline{z}}}{\lambda}\partial_z$, $<S\partial_z,\partial_z>=Q^{(2,0)}$ and $<S\partial_z,\partial_{\overline{z}}>=0$, where $\lambda=<\partial_z,\partial_{\overline{z}}>$.

Next, we compute
\begin{equation*}
<T_S(\partial_z,\partial_{\overline{z}}),\partial_{\overline{z}}>=\partial_{\overline{z}}<\partial_{\overline{z}},S\partial_z>-<S\partial_z,\nabla_{\partial_{\overline{z}}}\partial_{\overline{z}}>
\end{equation*}
\begin{equation*}
-\partial_z<S\partial_{\overline{z}},\partial_{\overline{z}}>+ <S\partial_{\overline{z}},\nabla_{\partial_z}\partial_{\overline{z}}>
\end{equation*}
\begin{equation*}
=-\overline{Q^{(2,0)}_{\overline{z}}}=0,
\end{equation*}
where we use that $\nabla_{\partial_{\overline{z}}}\partial_{\overline{z}}=\dfrac{\lambda_{\overline{z}}}{\lambda}\partial_{\overline{z}}$ and $\overline{ Q^{(2,0)}}_z=\overline{Q^{(2,0)}_{\overline{z}}}$.
It follows that $T_S=0$.

\cqd

\begin{lem}\label{l4}
Let $S$ be a symmetric operator satisfying Codazzi's equation and
$tr(S)=0$, then
\begin{equation}\label{e1}
<(\nabla^2S)x,y>=\sum_{i=1}^2\{-<Sy,R_{e_i,x}e_i>-<Se_i,R_{e_i,x}y>\},
\end{equation}
where $\{e_1,e_2\}$ is an orthonormal basis of \ $T_p\Sigma$.
\end{lem}
\dem \ \  Let us consider a geodesic frame $\{E_1,E_2\}$ in $p \in \Sigma^2$ which extend the basis $\{e_1,e_2\}$ and $X,Y$ local vector parallel fields which extend the given vectors $x,y$. Compute
\begin{equation*}
 (\nabla^2S)X= \sum_{i=1}^2(\nabla_{E_i}\nabla_{E_i}S)X= \sum_{i=1}^2\nabla_{E_i}((\nabla_{E_i}S)X)
\end{equation*}
\begin{equation*}
= \sum_{i=1}^2\nabla_{E_i}((\nabla_{X}S)E_i=\sum_{i=1}^2(\nabla_{E_i}\nabla_{X}S)E_i=
\sum_{i=1}^2(\nabla_{X}\nabla_{E_i}S)E_i+\sum_{i=1}^2(R(E_i,X)S)E_i
\end{equation*}
\begin{equation*}
=\sum_{i=1}^2\nabla_{X}((\nabla_{E_i}S)E_i)+\sum_{i=1}^2(R(E_i,X)S)E_i
\end{equation*}
\begin{equation*}
=\nabla_{X}\left(\sum_{i=1}^2(\nabla_{E_i}S)E_i\right)+\sum_{i=1}^2(R(E_i,X)S)E_i,
\end{equation*}
where we have used that $S$ satisfying Codazzi's  equation.

Then,
\begin{equation*}
 <(\nabla^2S)x,y>= <\nabla_{X}\left(\sum_{i=1}^2(\nabla_{E_i}S)E_i\right),Y>+\sum_{i=1}^2<(R(E_i,X)S)E_i,Y>
\end{equation*}
\begin{equation*}
= X\left(\sum_{i=1}^2<(\nabla_{E_i}S)E_i,Y>\right)+\sum_{i=1}^2<(R(E_i,X)S)E_i,Y>
\end{equation*}
\begin{equation*}
 =X(tr(\nabla_YS))+\sum_{i=1}^2<(R(E_i,X)S)E_i,Y>
\end{equation*}
\begin{equation*}
 =XY(tr(S))+\sum_{i=1}^2(<(R(E_i,X)(SE_i),Y>-<S(R(E_i,X)E_i),Y>).
\end{equation*}
By using that $tr S=0$, and computing the above at the point $p$, we obtain
\begin{equation*}
 <(\nabla^S)x,y>=\sum_{i=1}^2(-<(R(e_i,x)y,Se_i>-<R(e_i,x)e_i,Sy>).
\end{equation*}
 \cqd

Let us evaluate each summand in the expression $(\ref{e1})$.

\begin{lem}\label{l5}
Let $S$ be an operator as in Lemma $\ref{l4}$. Then,
\begin{equation*}i) \ \sum_{i=1}^2<Sy,R_{e_i,x}e_i>=-c\nu^2<Sx,y>-2H<Ax,Sy>+<A^2x,Sy>.\end{equation*}
and
\begin{equation*}ii) \ \sum_{i=1}^2<Se_i,R_{e_i,x}y>=-c\nu^2<Sx,y>-<Ay,SAx>+<Ax,y>tr(AS).\end{equation*}
\end{lem}
\dem \ \ Consider $\{e_1,e_2\}$ an orthonormal basis of
$T_p\Sigma$. Using the Gauss' equation $(\ref{gs})$ we have
\begin{equation*}
<Sy,R_{e_i,x}e_i>=-c\{<x,Sy>-<x,e_i><Sy,e_i>-<x,T><Sy,T>
\end{equation*}
\begin{equation*}
-<e_i,T>^2<x,Sy>+<e_i,T><x,e_i><Sy,T>
\end{equation*}
\begin{equation*}
+<x,T><e_i,T><e_i,Sy>\}-<Ae_i,e_i><Ax,Sy>
\end{equation*}
\begin{equation*}
+<Ax,e_i><Ae_i,Sy>.
\end{equation*}
Therefore,
\begin{equation*}
\sum_{i=1}^2<Sy,R_{e_i,x}e_i>=-c\{2<x,Sy>-\sum_{i=1}^2<x,e_i><Sy,e_i>
\end{equation*}
\begin{equation*}
-2<x,T><Sy,T>-<x,Sy>\sum_{i=1}^2<e_i,T>^2
\end{equation*}
\begin{equation*}
+<Sy,T>\sum_{i=1}^2<e_i,T><x,e_i>+<x,T>\sum_{i=1}^2<e_i,T><e_i,Sy>\}
\end{equation*}
\begin{equation*}
-<Ax,Sy>\sum_{i=1}^2<Ae_i,e_i>+\sum_{i=1}^2<Ax,e_i><Ae_i,Sy>,
\end{equation*}
which implies that
\begin{equation*}
\sum_{i=1}^2<Sy,R_{e_i,x}e_i>=-c\{2<x,Sy>-<Sx,y>-2<x,T><Sy,T>
\end{equation*}
\begin{equation*}
-<x,Sy>|T|^2+<Sy,T><x,T>+<x,T><T,Sy>\}
\end{equation*}
\begin{equation*}
-<Ax,Sy>2H+<Ax,ASy>.
\end{equation*}

Hence,
\begin{equation*}\sum_{i=1}^2<Sy,R_{e_i,x}e_i>=-c(1-|T|^2)<Sx,y>-2H<Ax,Sy>+<A^2x,Sy>,\end{equation*}
which shows the validity of $(i)$.
Now,
\begin{equation*}
<Se_i,R_{e_i,x}y>=-c\{<e_i,y><Se_i,x>-<x,y><Se_i,e_i>
\end{equation*}
\begin{equation*}
-<x,T><Se_i,T><e_i,y>-<e_i,T><y,T><x,Se_i>
\end{equation*}
\begin{equation*}
+<e_i,T><x,y><Se_i,T>+<x,T><y,T><e_i,Se_i>\}
\end{equation*}
\begin{equation*}
-<Ae_i,y><Ax,Se_i>+<Ax,y><Ae_i,Se_i>.
\end{equation*}
Therefore
\begin{equation*}
\sum_{i=1}^2<Se_i,R_{e_i,x}y>=-c\{\sum_{i=1}^2<e_i,y><Se_i,x>-<x,y>\sum_{i=1}^2<Se_i,e_i>
\end{equation*}
\begin{equation*}
-<x,T>\sum_{i=1}^2<Se_i,T><e_i,y>-<y,T>\sum_{i=1}^2<e_i,T><x,Se_i>
\end{equation*}
\begin{equation*}
+<x,y>\sum_{i=1}^2<e_i,T><Se_i,T>+<x,T><y,T>\sum_{i=1}^2<e_i,Se_i>\}
\end{equation*}
\begin{equation*}
-\sum_{i=1}^2<Ae_i,y><Ax,Se_i>+<Ax,y>\sum_{i=1}^2<Ae_i,Se_i>.
\end{equation*}
Therefore
\begin{equation*}\sum_{i=1}^2<Se_i,R_{e_i,x}y>=-c\{<Sx,y>-<x,T><Sy,T>-<y,T><Sx,T>\end{equation*}
\begin{equation*}+<ST,T><x,y>\}-<Ay,SAx>+<Ax,y>tr(AS),\end{equation*}
where we used that $trS=0$.

Using the fact that
\begin{equation*}-<x,T><Sy,T>-<y,T><Sx,T>+<ST,T><x,y>
\end{equation*}
\begin{equation*}=-(1-\nu^2)<Sx,y>,\end{equation*} we find that
\begin{equation*}\sum_{i=1}^2<Se_i,R_{e_i,x}y>=-c\nu^2<Sx,y>-<Ay,SAx>+<Ax,y>tr(AS),\end{equation*} which shows $(ii)$.
\cqd

\begin{cor}\label{c3}
Let $\Sigma^2 \looparrowright M^2(c)\times\real$ be an immersed
surface with constant mean curvature $H$ and $S$ the operator
satisfying Codazzi's equation with $trS=0$. Then,
\begin{equation*}
<(\nabla^2S)x,y>=2c\nu^2<Sx,y>+2H<Ax,Sy>-<A^2x,Sy>
\end{equation*}
\begin{equation*}
+<Ay,SAx>-<Ax,y>tr(AS).
\end{equation*}
\end{cor}
\dem \ \ We use the expressions of Lemma $\ref{l5}$ in
the equation $(\ref{e1})$ obtained in Lemma $\ref{l4}$. \cqd

Next we obtain an equation of Simons' type for the operator $S$ as
defined in $(\ref{e0})$.

\begin{thm}[Theorem 1.2 in Introduction]\label{t32}
Let $\Sigma^2 \looparrowright M^2(c)\times\real$ be an immersed
surface with non zero constant mean curvature $H$ and $S$ as
defined in $(\ref{e0})$. Then,
\begin{equation*}
\dfrac{1}{2}\Delta|S|^2=|\nabla
S|^2-|S|^4+|S|^2\left(\dfrac{5c\nu^2}{2}-\dfrac{c}{2}+2H^2-\dfrac{c}{H}<ST,T>\right)
\end{equation*}
\begin{equation*}
+c|ST|^2-\dfrac{1}{4H^2}<ST,T>^2.
\end{equation*}
\end{thm}
\dem \ \ We know that  $\dfrac{1}{2}\Delta|S|^2= |\nabla
S|^2+<\nabla^2S,S>$. Furthermore, using Corollary $\ref{c3}$, we
obtain
\begin{equation*}
<\nabla^2S,S>=2c\nu^2|S|^2+2Htr(AS^2)-[tr(AS)]^2.
\end{equation*}
Now, we need to compute $tr(AS^2)$ and $tr(AS)$.
\\
We have
\\
\begin{equation*}
tr(AS^2)=tr\{S^2(S+\dfrac{c}{2H}<T,\cdot>
T-\dfrac{c}{4H}(1-\nu^2)I+HI)\}
\end{equation*}
\begin{equation*}
=trS^3+\dfrac{c}{2H}tr(<T,S^2\cdot>
T)-\left(\dfrac{c}{4H}(1-\nu^2)-H\right)trS^2
\end{equation*}
\begin{equation*}
=0+\dfrac{c}{2H}|ST|^2-\left(\dfrac{c}{4H}(1-\nu^2)-H\right)|S|^2
\end{equation*}
and
\begin{equation*}
tr(AS)=tr\{S(S+\dfrac{c}{2H}<T,\cdot>
T-\dfrac{c}{4H}(1-\nu^2)I-HI)\}
\end{equation*}
\begin{equation*}
=trS^2+\dfrac{c}{2H}tr(<T,S\cdot> T)-(\dfrac{c}{4H}(1-\nu^2)-H)trS
\end{equation*}
\begin{equation*}
=|S|^2+\dfrac{c}{2H}<ST,T>-0,
\end{equation*}
where we used that $trS=trS^3=0$, that
\begin{equation*}
tr(<T,S\cdot> T)= \sum_{i=1}^2<T,Se_i><T,e_i>=<ST,T>
\end{equation*}
and that
\begin{equation*}
tr(<T,S^2\cdot> T)= \sum_{i=1}^2<T,S^2e_i><T,e_i>=<S^2T,T>.
\end{equation*}

 Therefore,
\begin{equation*}
\dfrac{1}{2}\Delta|S|^2= |\nabla
S|^2+2c\nu^2|S|^2+2H\left(\dfrac{c}{2H}|ST|^2-\left(\dfrac{c}{4H}(1-\nu^2)-H\right)|S|^2\right)
\end{equation*}
\begin{equation*}
-\left(|S|^2+\dfrac{c}{2H}<ST,T>\right)^2,
\end{equation*}
so that,
\begin{equation*}
\dfrac{1}{2}\Delta|S|^2= |\nabla
S|^2+2c\nu^2|S|^2+c|ST|^2-\left(\dfrac{c}{2}(1-\nu^2)-2H^2\right)|S|^2
\end{equation*}
\begin{equation*}
-|S|^4-\dfrac{c}{H}<ST,T>|S|^2-\dfrac{1}{4H^2}<ST,T>^2.
\end{equation*}

Rearranging terms, we obtain finally
\begin{equation*}
\dfrac{1}{2}\Delta|S|^2=|\nabla
S|^2-|S|^4+|S|^2\left(\dfrac{5c\nu^2}{2}-\dfrac{c}{2}+2H^2-\dfrac{c}{H}<ST,T>\right)
\end{equation*}
\begin{equation*}
+c|ST|^2-\dfrac{1}{4H^2}<ST,T>^2.
\end{equation*}
\cqd


\section{Applications}

In this section, we will apply the results in section 3 together
with the Omori-Yau's Theorem and classify some surfaces in
$M^2(c)\times\real$.

\begin{thm}
Let $\Sigma^2 \looparrowright \h\times\real$ be an oriented
complete non-compact immersed minimal surface. Assume that
\begin{equation*}
\sup_{\Sigma}(|A|^2+5\nu^2) < 1.
\end{equation*}
Then $\Sigma^2$ is a vertical plane $\gamma\times\real$ for some
geodesic $\gamma$ in $\h$.
\end{thm}
\dem \ \  Using Theorem $\ref{t31}$ with $H=0$ and
$c=-1$, so
\begin{equation*}
\dfrac{1}{2}\Delta|A|^2= |\nabla A|^2 -|A|^4+(1-5\nu^2)|A|^2\geq
|A|^2\left(-|A|^2+1-5\nu^2\right).
\end{equation*}
Let $\dfrac{d}{2}:=  -\sup_{\Sigma}(|A|^2+5\nu^2) + 1>0$.
Therefore,
\begin{equation}\label{e2}
\Delta|A|^2 \geq d\cdot|A|^2.
\end{equation}
Using Gauss' equation $(\ref{gs})$ in $\h\times\real$ we have
\begin{equation*}
K_{\Sigma}=K_{ext}-\nu^2 =
-\dfrac{|A|^2+5\nu^2}{2}+\dfrac{3\nu^2}{2}\geq-\dfrac{1}{2}.
\end{equation*}
Now we can use Theorem $\ref{OY}$ with $u=|A|^2$, i.e, there exist
$\{p_j\}$ in $\Sigma^2$ such that
\begin{equation*}
\lim_{j\rightarrow\infty}|A|^2(p_j)=\sup_{\Sigma} |A|^2 \ and \
\lim_{j\rightarrow\infty}\Delta |A|^2(p_j) \leq 0.
\end{equation*}
Therefore we use inequality $(\ref{e2})$ to conclude that
$\sup_{\Sigma}|A|^2=0$, i.e, $\Sigma^2$ is totally geodesic with
$|\nu| < \sqrt{0.2}$.
\newpage
Since $\Sigma^2$ is totally geodesic and $|\nu| < \sqrt{0.2}$ it
cannot be a slice, it must be a vertical plane $\gamma\times\real$
for some geodesic $\gamma$ in $\h$.

 \cqd


\begin{thm}
Let $\Sigma^2 \looparrowright \h\times\real$ be a complete
immersed surface with constant mean curvature $H$. Assume that
\begin{equation*}
\sup_{\Sigma}(| \phi|^2+5 \nu^2)<2H^2+1 \ \ and \ \  <\phi T,T>
\geq 0.
\end{equation*}
Then $\Sigma^2$ is a vertical plane $\gamma\times\real$ for some
geodesic $\gamma$ in $\h$.
\end{thm}
\dem \ \ We consider the expression in Theorem $\ref{t31}$ with
$c=-1$:
\begin{equation*}
\dfrac{1}{2}\Delta|\phi|^2= |\nabla\phi|^2
-|\phi|^4+(2H^2+1-5\nu^2)|\phi|^2+2H<\phi T,T>.
\end{equation*}
As $<\phi T,T>\geq 0$, we have
\begin{equation*}
\dfrac{1}{2}\Delta|\phi|^2 \geq -|\phi|^4+(2H^2+1-5\nu^2)|\phi|^2.
\end{equation*}
Consider $\dfrac{d}{2}:= 2H^2+1-\sup_{\Sigma}(| \phi|^2+5
\nu^2)>0$. Then
\begin{equation*}
\Delta|\phi|^2 \geq 2|\phi|^2(2H^2+1-5\nu^2-|\phi|^2)\geq
d|\phi|^2,
\end{equation*}
which implies,
\begin{equation}\label{e3}
\Delta|\phi|^2 \geq d|\phi|^2.
\end{equation}

Using Gauss' equation $(\ref{gs})$ in $\h\times\real$ we have
\begin{equation*}
K_{\Sigma}=K_{ext}-\nu^2 =
-\dfrac{|\phi|^2+5\nu^2-2H^2}{2}+\dfrac{3\nu^2}{2}\geq-\dfrac{1}{2}.
\end{equation*}
Now we can use Theorem $\ref{OY}$ with $u=|\phi|^2$, i.e, there
exist $\{p_j\}$ in $\Sigma^2$ such that
\begin{equation*}
\lim_{j\rightarrow\infty}|\phi|^2(p_j)=\sup_{\Sigma} |\phi|^2 \
and \ \lim_{j\rightarrow\infty}\Delta |\phi|^2(p_j) \leq 0.
\end{equation*}
Furthermore, we use inequality $(\ref{e3})$ to conclude that
$\sup_{\Sigma}|\phi|^2=0$, i.e, $\Sigma^2$ is totally umbilical.
Next, we use the results obtained in $[ST]$ in section 4 to conclude
that $\Sigma^2$ is totally geodesic, where we have used that if
$\Sigma^2$ is totally umbilical with constant mean curvature in
$\h\times\real$ then $\Sigma^2$ is totally geodesic.
\newpage
Since $\Sigma^2$ is totally geodesic and $|\nu| < \sqrt{0.2}$ it
must be a vertical plane $\gamma\times\real$ for some geodesic
$\gamma$ in $\h$. This concludes the proof.\cqd


 We need the following result:

\begin{lem}\label{l41}
Let $\Sigma^2\looparrowright M^2(c)\times\real$ be a complete
immersed surface with non zero constant mean curvature $H$. Then
$|S|=0$ if and only if $\Sigma^2$ is an Abresch-Rosenberg surface.
\end{lem}
\dem \ \ We consider $(u,v)$ an isothermal parameter to a surface
$\Sigma^2$. Now, we consider the complex parameter to
riemannian metric, $z=u+iv$ and (2,0)-part of the Abresch-Rosenberg
differential
\begin{equation*}Q(u,v)=2H<Au,v>-c<u,T><v,T>.\end{equation*}
 We can rewrite $Q$ as
\begin{equation*}
Q(u,v)=<Su,v>-\dfrac{c}{2}(1-\nu^2)<u,v>+2H^2<u,v>.
\end{equation*}
Next we evaluate $Q(\partial_z,\partial_z)$ and use that
$<\partial_z,\partial_z>=0$:
\begin{equation*}
Q(\partial_z,\partial_z)=<S\partial_z,\partial_z>=\left(\dfrac{\tilde{e}-\tilde{g}}{4}\right)-i\dfrac{\tilde{f}}{2},
\end{equation*}
where
$\tilde{e}=<S\partial_u,\partial_u>=-<S\partial_v,\partial_v>=-\tilde{g}$
and $\tilde{f}=<S\partial_u,\partial_v>$.  Therefore

\begin{equation*}
|Q^{(2,0)}|=\sqrt{\left(\dfrac{\tilde{e}-\tilde{g}}{4}\right)^2+\dfrac{\tilde{f}^2}{4}}=\sqrt{\dfrac{\tilde{e}^2}{4}+\dfrac{\tilde{f}^2}{4}}=\dfrac{1}{2\sqrt{2}}|S|.
\end{equation*}
This concludes the proof. \cqd


Let us consider the polynomial
$p_H(t)=-t^2-\dfrac{1}{H}t+\left(\dfrac{4H^2-1}{2}\right)$.  When
$H$ is greater than one half then there is a positive root for
$p_H$. Let $L_H$ be the positive root. One has:

\begin{thm}[Theorem 1.3 in Introduction]
Let $\Sigma^2 \looparrowright \s\times\real$ be a complete
immersed surface with constant mean curvature $H$ greater than one
half. Assume that
\begin{equation*}
\sup_{\Sigma}|S| < L_H.
\end{equation*}
Then $\Sigma^2 = S_H^2$, i.e, $\Sigma^2$ is a Hsiang-Pedrosa
sphere.
\end{thm}
\dem \ \ Consider the expression in Theorem $\ref{t32}$ with
$c=1$:
\begin{equation*}
\dfrac{1}{2}\Delta|S|^2=|\nabla
S|^2-|S|^4+|S|^2\left(\dfrac{5\nu^2}{2}-\dfrac{1}{2}+2H^2-\dfrac{1}{H}<ST,T>\right)
\end{equation*}
\begin{equation*}
+|ST|^2-\dfrac{1}{4H^2}<ST,T>^2.
\end{equation*}
As $|<ST,T>|\leq |ST|\leq |S|$, we have
\begin{equation*}
\dfrac{1}{2}\Delta|S|^2\geq-|S|^4+|S|^2\left(\dfrac{4H^2-1}{2}-\dfrac{1}{H}|S|\right)+\left(\dfrac{4H^2-1}{4H^2}\right)<ST,T>^2,
\end{equation*}
hence,
\begin{equation}\label{n1}
\dfrac{1}{2}\Delta|S|^2\geq
|S|^2\left(\dfrac{4H^2-1}{2}-\dfrac{1}{H}|S|-|S|^2\right)+\dfrac{5}{2}\nu^2|S|^2,
\end{equation}
because $H>\dfrac{1}{2}$.
\\
Observe that
\begin{equation*}
\dfrac{4H^2-1}{2}-\dfrac{1}{H}|S|-|S|^2 \geq
p_H(sup_{\Sigma}|S|)=:\dfrac{d}{2}>0
\end{equation*}
and $\nu^2|S|^2\geq 0$. Therefore
\begin{equation}\label{e4}
\Delta|S|^2\geq d|S|^2.
\end{equation}
Now we estimate $|S|$.
\begin{equation*}
|S|\geq2H|A|-|<T,\cdot>T|- (1-\nu^2)-4H^2\geq
2H|A|-2(1-\nu^2)-4H^2,
\end{equation*}
that is,
\begin{equation*}
L_H\geq |S|\geq2H|A|- 2-4H^2.
\end{equation*}
Using Gauss' equation $(\ref{gs})$ in $\s\times\real$ we have
\begin{equation*}
K_{\Sigma}=K_{ext}+\nu^2 = -\dfrac{|A|^2}{2}+2H^2+\nu^2\geq
-\dfrac{1}{2}\left(\dfrac{L_H+2+4H^2}{2H}\right)^2.
\end{equation*}
Now we can use Theorem $\ref{OY}$ with $u=|S|^2$, i.e, there
exists a $\{p_j\}$ in $\Sigma^2$ such that
\begin{equation*}
\lim_{j\rightarrow\infty}|S|^2(p_j)=\sup_{\Sigma} |S|^2 \ and \
\lim_{j\rightarrow\infty}\Delta |S|^2(p_j) \leq 0.
\end{equation*}
We then use inequality $(\ref{e4})$ to conclude that
$\sup_{\Sigma}|S|^2=0$, i.e, $|S|=0$ in $\Sigma^2$. Using Lemma
$\ref{l41}$ we conclude the proof. \cqd

\begin{thm}
There exists no $\Sigma^2\looparrowright\s\times\real$ complete
immersed surface with constant mean curvature greater than one half
such that $|S|=L_H$.
\end{thm}
\dem \ \ Suppose that there exist
$\Sigma^2\looparrowright\s\times\real$ satisfying condition of the
theorem. Using the expression $(\ref{n1})$
\begin{equation*}
\dfrac{1}{2}\Delta|S|^2\geq
|S|^2\left(\dfrac{4H^2-1}{2}-\dfrac{1}{H}|S|-|S|^2\right)+\dfrac{5}{2}\nu^2|S|^2,
\end{equation*}
with $|S|=L_H$ we have:
\begin{equation*}
0\geq 0+\dfrac{5}{2}\nu^2L_H^2\geq 0.
\end{equation*}
Hence $\nu=0$, i.e, $\Sigma^2\looparrowright\s\times\real$ is a
cylinder $\gamma\times\real$ for some $\gamma \in \s$ with
constant curvature $2H$.

On the other hand, for the cylinder $\gamma\times\real$, where
$\gamma \in \s$ with constant curvature $2H$, we have that
\begin{equation*}
S=\left(%
\begin{array}{cc}
  2H^2+\dfrac{1}{2} & 0 \\
  0 & -2H^2-\dfrac{1}{2} \\
\end{array}%
\right).
\end{equation*}
As $|S|=\dfrac{\sqrt{2}}{2}(4H^2+1) > L_H$ we have a
contradiction. \cqd

\begin{thm}
Let $\Sigma^2 \looparrowright \s\times\real$ be a closed immersed
surface with constant mean curvature $H$ greater than one half.
Assume that
\begin{equation*}
|S| \leq  L_H.
\end{equation*}
Then $\Sigma^2 = S_H^2$, i.e, $\Sigma^2$ is a Hsiang-Pedrosa
sphere.
\end{thm}
\dem \ \ Let us consider the expression $(\ref{n1})$
\begin{equation*}
\dfrac{1}{2}\Delta|S|^2\geq
|S|^2\left(\dfrac{4H^2-1}{2}-\dfrac{1}{H}|S|-|S|^2\right)+\dfrac{5}{2}\nu^2|S|^2.
\end{equation*}
As $|S| \leq  L_H$, we have
$\dfrac{4H^2-1}{2}-\dfrac{1}{H}|S|-|S|^2 \geq 0$. Hence,
\begin{equation*}
\dfrac{1}{2}\Delta|S|^2\geq \dfrac{5}{2}\nu^2|S|^2.
\end{equation*}
Integrating and using Stokes' Theorem we have
\begin{equation*}
0\geq \dfrac{5}{2}\int_{\Sigma}\nu^2|S|^2 d\Sigma\geq 0.
\end{equation*}
It follows that $|S|\cdot\nu=0$. If $\nu=0$ then we have a cylinder,
but this is not possible because $\Sigma^2$ is closed. Therefore
$|S|=0$. Using Lemma $\ref{l41}$ we conclude the proof.

\cqd

\vspace{1cm}

 Let us consider the polynomial
$q_H(t)=-t^2-\dfrac{1}{H}t+\left(\dfrac{8H^4-12H^2-1}{4H^2}\right)$.
When $H$ is greater than a positive root of the polynomial
$r(x)=8x^4-12x^2-1$, i.e, $H$ greater than
$\sqrt{\dfrac{12+\sqrt{176}}{16}}$, then there is a positive root
for $q_H$. Let $M_H$ be the positive root.
\\
\begin{thm}[Theorem 1.4 in Introduction]
Let $\Sigma^2 \looparrowright \h\times\real$ be a complete
immersed surface with constant mean curvature $H$ greater than
$\sqrt{\dfrac{12+\sqrt{176}}{16}}\approx 1.25664$. Assume that
\begin{equation*}
\sup_{\Sigma}|S| < M_H.
\end{equation*}
Then  $\Sigma^2$ is an Abresch-Rosenberg surface.
\end{thm}
\dem \ \ Consider the expression in Theorem $\ref{t32}$ with
$c=-1$
\begin{equation*}
\dfrac{1}{2}\Delta|S|^2=|\nabla
S|^2-|S|^4+|S|^2\left(-\dfrac{5\nu^2}{2}+\dfrac{1}{2}+2H^2+\dfrac{1}{H}<ST,T>\right)
\end{equation*}
\begin{equation*}
-|ST|^2-\dfrac{1}{4H^2}<ST,T>^2.
\end{equation*}
As $|<ST,T>|\leq |ST|\leq |S|$, we have
\begin{equation*}
\dfrac{1}{2}\Delta|S|^2\geq-|S|^4+|S|^2\left(\dfrac{4H^2+1-5\nu^2}{2}-\dfrac{1}{H}|S|\right)-\left(\dfrac{4H^2+1}{4H^2}\right)|S|^2,
\end{equation*}
i.e,
\begin{equation*}
\dfrac{1}{2}\Delta|S|^2\geq
|S|^2\left(\dfrac{4H^2-4+5-5\nu^2}{2}-\dfrac{1}{H}|S|-\dfrac{4H^2+1}{4H^2}-|S|^2\right).
\end{equation*}
This may be rewritten as,
\begin{equation}\label{n2}
\dfrac{1}{2}\Delta|S|^2\geq
|S|^2\left(\dfrac{8H^4-12H^2-1}{4H^2}-\dfrac{1}{H}|S|-|S|^2\right)+\dfrac{5}{2}(1-\nu^2)|S|^2.
\end{equation}
\\
Observe that
\begin{equation*}
\dfrac{8H^4-12H^2-1}{4H^2}-\dfrac{1}{H}|S|-|S|^2 \geq
q_H(sup_{\Sigma}|S|)=:\dfrac{d}{2}>0
\end{equation*}
 and $(1-\nu^2)|S|^2\geq 0$. Therefore,
\begin{equation}\label{e5}
\Delta|S|^2\geq d|S|^2.
\end{equation}
Next we estimate $|S|$.
\begin{equation*}
|S|\geq2H|A|-|<T,\cdot>T|- (1-\nu^2)-4H^2\geq
2H|A|-2(1-\nu^2)-4H^2,
\end{equation*}
i.e,
\begin{equation*}
M_H\geq |S|\geq2H|A|- 2-4H^2.
\end{equation*}
Using Gauss' equation $(\ref{gs})$ in $\h\times\real$ we have
\begin{equation*}
K_{\Sigma}=K_{ext}-\nu^2 = -\dfrac{|A|^2}{2}+2H^2-\nu^2\geq
-\dfrac{1}{2}\left(\dfrac{M_H+2+4H^2}{2H}\right)^2.
\end{equation*}
Now we can use Theorem $\ref{OY}$ with $u=|S|^2$, i.e, there
exists a $\{p_j\}$ in $\Sigma^2$ such that
\begin{equation*}
\lim_{j\rightarrow\infty}|S|^2(p_j)=\sup_{\Sigma} |S|^2 \ and \
\lim_{j\rightarrow\infty}\Delta |S|^2(p_j) \leq 0.
\end{equation*}
We use inequality $(\ref{e5})$ to conclude that
$\sup_{\Sigma}|S|^2=0$, i.e, $|S|=0$ in $\Sigma^2$. Then using Lemma
$\ref{l41}$ we conclude the proof. \cqd \vspace{1cm}

\begin{thm}
There exists no $\Sigma^2\looparrowright\h\times\real$ a complete
immersed surface with constant mean curvature greater than
$\sqrt{\dfrac{12+\sqrt{176}}{16}}\approx 1.25664$  such that
$|S|=M_H$.
\end{thm}
\dem \ \ Suppose that there exists
$\Sigma^2\looparrowright\h\times\real$ satisfying condition of the
theorem. Using the expression $(\ref{n2})$
\begin{equation*}
\dfrac{1}{2}\Delta|S|^2\geq
|S|^2\left(\dfrac{8H^4-12H^2-1}{4H^2}-\dfrac{1}{H}|S|-|S|^2\right)+\dfrac{5}{2}(1-\nu^2)|S|^2
\end{equation*}
with $|S|=M_H$ we have:
\begin{equation*}
0\geq 0+\dfrac{5}{2}(1-\nu^2)M_H^2\geq 0.
\end{equation*}
Hence $\nu^2=1$, i.e, $\Sigma^2\looparrowright\h\times\real$ is a
slice $\h\times\{t_0\}$. But $\h\times\{t_0\}$ has zero mean
curvature, and this is impossible because $H$ is positive.

\cqd

\begin{thm}
Let $\Sigma^2 \looparrowright \h\times\real$ be a closed immersed
surface with constant mean curvature $H$ greater than
$\sqrt{\dfrac{12+\sqrt{176}}{16}}\approx 1.25664$ . Assume that
\begin{equation*}
|S| \leq  M_H.
\end{equation*}
Then $\Sigma^2 = S_H^2$, i.e, $\Sigma^2$ is a Hsiang-Pedrosa
sphere.
\end{thm}
\dem \ \ Let us consider the expression $(\ref{n2})$
\begin{equation*}
\dfrac{1}{2}\Delta|S|^2\geq
|S|^2\left(\dfrac{8H^4-12H^2-1}{4H^2}-\dfrac{1}{H}|S|-|S|^2\right)+\dfrac{5}{2}(1-\nu^2)|S|^2.
\end{equation*}
As $|S| \leq  M_H$, we have
$\dfrac{8H^4-12H^2-1}{4H^2}-\dfrac{1}{H}|S|-|S|^2 \geq 0$. Hence,
\begin{equation*}
\dfrac{1}{2}\Delta|S|^2\geq \dfrac{5}{2}(1-\nu^2)|S|^2.
\end{equation*}
Integrating and using Stokes' Theorem we have
\begin{equation*}
0\geq \dfrac{5}{2}\int_{\Sigma}(1-\nu^2)|S|^2 d\Sigma\geq 0.
\end{equation*}
Moreover $(1-\nu^2)\cdot|S|^2=0$. If $\nu^2=1$, again we have a
slice, but this is not possible. Therefore $|S|=0$. Using Lemma
$\ref{l41}$ we conclude the proof.

\cqd


\hspace{-0.6cm} Márcio Batista \\
Universidade Federal de Alagoas\\
Instituto de Matemática\\
57072-900 -  Maceió - Alagoas - Brazil\\
$\it{email-address}$: mhbs28@gmail.com

\end{document}